\documentclass[a4paper,12pt]{article}
\usepackage{amsfonts,amsmath,amssymb,amscd,amsthm,latexsym}

\let\cal\mathcal
\newcommand{\bi}{\bibitem}
\newcommand{\nb}{\newblock}
\newcommand{\be}[1]{\begin{equation}\label{#1}}
\newcommand{\ee}{\end{equation}}
\newcommand{\la}{\langle\,}
\newcommand{\ra}{\,\rangle}
\newcommand{\ve}{\varepsilon}
\newcommand{\ccc}{{\cal C}}
\newcommand{\Card}{\mathop{\rm Card}}

\newtheorem{thm}{\quad Theorem}
\newtheorem{lm}{\quad Lemma}

\title{Traveller Salesman Property and Richard Thompson's Group $F$}
\author{\vspace{2ex}
V. S. Guba\thanks{This research is partially supported by the RFFI
grant 05-01-00895.}\\
Vologda State Pedagogical University,\\
6 S.\,Orlov Street,\\
Vologda\\
Russia\\
160600\\
E-mail: guba{@}uni-vologda.ac.ru}
\date{}

\begin{document}

\maketitle

%\subjclass{20F32,05C25}
%\keywords{Richard Thompson's group $F$;amenability;traveller salesman
%property}

\begin{abstract}

Recently Akhmedov introduced a property of finitely generated
groups called the Traveller Salesman Property (TSP). It was shown
by Thurston that TSP implies non-amenability. Akhmedov conjectured
that R.\,Thompson's group $F$ has TSP. The aim of this article is
to disprove this conjecture.

\end{abstract}

The Traveller Salesman Property of groups was introduced by
Akhmedov in his PhD Thesis; see also \cite{Akh05}. We start with
several definitions.

Let $G$ be a group and let $\xi\in G$. A subset $S\subseteq G$ is
called a $\xi$-{\em related\/} set whenever for any $g\in S$ one
has $g\xi\in S$ or $g\xi^{-1}\in S$. Elements of the form $g$
and $g\xi^{\pm1}$ are called $\xi$-{\em neighbours\/}.

Let $\Gamma$ be any connected graph and let $S$ be any finite
nonempty set of its vertices. Let $p$ be a shortest closed path
in $\Gamma$ that visits each vertex of $S$ at least once. By
$\tau(S)$ we denote the ratio $|p|/\Card(S)$, where $|p|$ is
the length of $p$.

Let $G$ be a group generated by a finite set $A$. By $\ccc(G,A)$
we denote the (right) Cayley graph of $G$ with $A$ as a set
of generators. For any $\lambda>0$, the graph $\ccc(G,A)$ is
said to satisfy property $TS(\lambda)$ whenever there exists a
$\xi\in G$ such that for any nonempty finite $\xi$-related set
$S$ one has $\tau(S)\ge\lambda$.

If $\ccc(G,A)$ satisfies $TS(\lambda)$ for all $\lambda>0$, then
we say that $\ccc(G,A)$ satisfies $TS(\infty)$. Notice that the
property $TS(\lambda)$ depends on the generating set of $G$.
However, it is easy to see that the property $TS(\infty)$ does
not depend on the choice of a generating set of $G$. Hence one
can say that a finitely generated group $G$ {\em has property\/}
$TS$ if $\ccc(G,A)$ satisfies $TS(\infty)$ for some (and so for
any) finite generating set $A$.

It was mentioned by Thurston that property $TS$ easily implies
non-amenability. For the sake of completeness, let us give an
argument that is close to Thurston's.

Suppose that a finitely generated group $G$ is amenable.
Let $A$ be any finite generating set of $G$. Given a subset
$S$ in $G$, by $\partial S$ we denote the boundary of $S$,
that is, the set of edges connecting a vertex in $S$ to a
vertex in $G\setminus S$. If $S$ is finite and nonempty,
then $i_*(S)$ denotes the isoperimetric constant of $S$, that
is, $i_*(S)=\Card(\partial S)/\Card(S)$. Using the F\o{}lner
criterion \cite{Fol}, for any $\ve>0$ one can find a set
$S=S_\ve$ satisfying $i_*(S)<\ve$. Notice that to every set
$S\subseteq G$ one can assign the corresponding subgraph in
$\ccc(G,A)$. In our case, without loss of generality one can
assume that the subgraph assigned to $S$ is connected (it
suffices to take a connected component with the smallest
value of $i_*$). So one can find a spanning subtree $T$ in
the subgraph formed by $S$. Thus there exists a closed path
of length $2\Card(S)$ visiting all vertices of $S$.

Given a $\xi\in G$, the number of vertices of $S$ that do
not have $\xi$-neighbours in $S$, does not exceed $C\ve\Card(S)$,
where $C$ is a constant that depends on the length of $\xi$
and the cardinality of $A$ but does not depend on $\ve$.
Removing these elements from $S$ gives us a $\xi$-related
set $S'$ of cardinality at least $(1-C\ve)\Card(S)$. So we
have $\tau(S')\le2/(1-C\ve)$. Since $\ve=\ve(\xi)$ can be chosen
as small as we want, this means that $\ccc(G,A)$ does not
satisfy $TS(\lambda)$ provided $\lambda>2$. In particular,
$G$ does not have $TS$.
\vspace{1ex}

In several conference talks, Akhmedov conjectured that Richard
Thompson's group $F$ has property $TS$. If this were true, this
would imply non-amenability of $F$. Recall that $F$ has the
following group presentation:

\be{xinf}
\la x_0,x_1,x_2,\ldots\mid x_j{x_i}=x_ix_{j+1}\ (i<j)\,\ra.
\ee
This group can be generated by $x_0$, $x_1$. In terms of these
generators, it has a finite presentation of the form

\be{xfin}
\la x_0,x_1\mid x_1^{x_0^2}=x_1^{x_0x_1},x_1^{x_0^3}=x_1^{x_0^2x_1}\ra,
\ee
where $a^b=b^{-1}ab$ by definition.
\vspace{1ex}

It was shown by Brin and Squier \cite{BrSq} that $F$ has no free
non-abelian subgroups. On the other hand, it is known that $F$ is
not elementary amenable (that is, it cannot be constructed from
finite and abelian groups using certain natural group theoretical
operations that preserve amenability). The reader can find many
details about this group in the survey \cite{CFP}; for additional
information see \cite{BG,GbS,GuSa99,Gu04,Gu05}.

The question about amenability of $F$ was asked by Ross Geoghegan
in the 70s. It is still open (for the definition of amenability of
groups and its basic properties see \cite{GrL}). We are going to
show that $F$ does not have property $TS$. So the fact that $TS$
implies non-amenabilty does not help to solve the above open problem.
\vspace{1ex}

An equivalent definition of $F$ can be given as follows. Let us
consider all strictly increasing continuous piecewise-linear
functions from the closed unit interval onto itself. Take only those
of them that are differentiable except at finitely many dyadic
rational points and such that all slopes (derivatives) are integer
powers of $2$. These functions form a group under composition. This
group is isomorphic to $F$. (An explicit form of this representation
of $F$ can be found in \cite{CFP}.)

We begin with an elementary lemma which seems to be interesting
in itself. Let $u=x_1$, $v=x_0x_1x_0^{-2}$ be the elements of $F$.
Using the representation of $F$ by piecewise-linear functions on $[0;1]$,
one can easily check that the restriction of $u$ on $[0;1/2]$ is the
identity function, as well as the restriction of $v$ on $[1/2;1]$. In
particular, $u$ and $v$ commute.

\begin{lm}
\label{mixed}
Let $\xi$ be any element of F. Then there exists $\ve=\pm1$ such that the
following relation holds: $(\xi^\ve u\xi^{-\ve})v=v(\xi^\ve u\xi^{-\ve})$,
where $u=x_1$, $v=x_0x_1x_0^{-2}$. That is, $\xi^\ve u\xi^{-\ve}$ and $v$
commute in $F$.
\end{lm}

\proof
All elements of $F$ act on the right as functions from $[0;1]$ onto
itself. If $(1/2)\xi>1/2$, then we choose $\ve=-1$, otherwise $\ve=1$.
Obviously, $(1/2)\xi^\ve$ does not exceed $1/2$. For simplicity, we will
put $z=\xi^\ve$.

We know that $z$ takes $[0;1/2]$ into itself. Therefore, $zuz^{-1}$
restricted on $[0;1/2]$ is the identity function. Indeed, if
$t\in[0;1/2]$, then $tz$ also belongs to $[0;1/2]$. The function
$u=x_1$ acts on $tz$ identically, that is, $tzu=tz$. Hence
$tzuz^{-1}=t$.

The function $v$ is the identity function on $[1/2,1]$. Since $zuz^{-1}$
and $v$ have disjoint supports, they commute. This completes the proof.
\endproof
\vspace{2ex}

One can extract as a corollary that the group $F$ satisfies
the following ``mixed identity":

\be{mxdid}
[[gx_1g^{-1},x_0x_1x_0^{-2}],[g^{-1}x_1g,x_0x_1x_0^{-2}]]=1 \ee for
any $g\in F$. Here $[a,b]=a^{-1}a^b=a^{-1}b^{-1}ab$ by definition.
Notice that not every group satisfies nontrivial mixed identities.
(For instance, a free non-abelian group does not.) Notice also that
$F$ does not satisfy any group law \cite{BrSq}.

\begin{thm}
\label{fnotints}
Let $\ccc$ be the Cayley graph of $F$ generated by $x_0$, $x_1$. Suppose
that $\lambda>2.5$. Then for any $\xi\in F$ there exists a $\xi$-related
finite set $S$ such that there is a closed path of length strictly
less than $\lambda\Card(S)$ in $\ccc$ that visits all vertices of
$S$. In particular, $F$ does not have property $TS$.
\end{thm}

\proof
Let $\xi\in F$. Using Lemma \ref{mixed}, we take $z=\xi^{\pm1}$ such that
$zuz^{-1}$ commutes with $v$. Thus $u$ commutes with $z^{-1}vz$. Clearly,
$u=x_1$ has length $1$ and $v=x_0x_1x_0^{-2}$ has length $4$.

We fix any even natural number $n$ and define a finite set $S_n$ of
vertices in $\ccc$. This set will be $\xi$-related and it will have $2N^2$
elements, where $N=n+1$.

First let us define the set $S_n$ formally. Then we will show its
geometric structure to work with it conveniently. For simplicity, let
$w=zuz^{-1}$.

We take all elements of the form $v^pw^q$, where $0\le p,q\le n$. They
form a set $S_n'$ of cardinality $N^2$ (all these elements are different
because $v$, $w$ have disjoint supports). By definition, $S_n''=S_n'z$.
The set $S_n$ is the union of $S_n'$ and $S_n''$.

The union of $S_n'$ and $S_n''$ is almost always disjoint. An
exception is only the case when $z$ belongs to the (free abelian)
subgroup generated by $v$ and $w$. This implies that $z$ commutes
with $w=zuz^{-1}$ and so it commutes with $u$ also. Hence $w=u$. Now
we know that $z$ is an element of the free abelian subgroup
generated by $u$ and $v$, the elements of bounded length. In this
case to construct a $\xi$-related set with desired properties is
very easy. Indeed, one can take a set $\bar S_n$ of all elements
of the form $u^pv^q$, where $0\le p,q\le n$. If $n$ is odd, then
the path labelled by $u^nvu^{1-n}vu^{n-1}v\cdots vu^{n-1}vu^{-n}v^{-n}$
visits all vertices of $\bar S_n$ and has length $n^2+8n+1$. If
$n$ is big enough, then the amount of elements that have no
$\xi$-neighbours in $\bar S_n$, is very small with respect to
$\Card(\bar S_n)=(n+1)^2$. This implies that for any $\lambda>1$
there exists a $\xi$-related set $\tilde S_n$ satisfying
$\tau(\tilde S_n)<\lambda$.
\vspace{1ex}

From now on, we can assume that the set $S_n=S_n'\cup S_n''$ has
cardinality $2N^2$. It is $\xi$-related because the map $g \mapsto gz$
is a bijection from $S_n'$ onto $S_n''$ such that $g$ and $gz$ are
$\xi$-neighbours.

Before constructing the closed path that visits all vertices of $S_n$, let
us try to imagine this set geometrically in the $3$-dimensional coordinate
space. Let us draw $n+1$ paths $p_0$, \dots, $p_n$ on the coordinate plane
with zero 3rd coordinate. Here $p_i$ starts at $(0,i,0)$ and ends at
$(n,i,0)$, where $0\le i\le n$. Every path $p_i$ is a product of $n$ edges,
each labelled by $v$. Now we draw $n+1$ paths $q_0$, \dots, $q_n$ on the
coordinate plane with the 3rd coordinate equal to $1$. For every
$0\le i\le n$, the endpoints of $q_i$ are $(i,0,1)$ and $(i,n,1)$.
Similarly, every $q_i$ is a product of $n$ edges, each labelled by $u$.

Finally, for every pair $(i,j)$ we put an arrow $r_{i,j}$ from $(i,j,0)$ to
$(i,j,1)$ and give the label $z$ to it (and also subdivide it into $|z|$
edges). The labelled graph thus obtained will be denoted by $\Gamma_n$.

We advise the reader to draw the graph $\Gamma_n$ for the case, say,
$n=4$. Clearly, we have a natural graph morphism from $\Gamma_n$ to
the Cayley graph $\ccc$. This morphism preserves all labels; it sends
$(i,j,0)$ to $v^iw^j$ and $(i,j,1)$ to $v^iw^jz$, which is also equal
to $z(w')^iu^j$, where $w'=z^{-1}vz$ by definition. Indeed,
$v^iw^jz=zz^{-1}v^iw^jz=z(z^{-1}vz)^i(z^{-1}wz)^j=z(w')^iu^j$.
Recall that $v$ commutes with $w$ and so $u$ commutes with $w'$.

Now we are ready to define a closed path $r$ in $\Gamma_n$ that will visit
each of the vertices of the form $(i,j,0)$ and $(i,j,1)$. Its image under
the graph morphism will have the same length $|r|$ and it will visit each
vertex of $S_n$. Here is the description:

\begin{itemize}
\item
we start at (0,0,0);
\item
go along $p_0$;
\item
go up by the arrow labelled by $z$;
\item
go along the first edge of $q_n$;
\item
go down by $z^{-1}$ (now we are at the point $(n,1,0)$);
\item
go along $p_1^{-1}$;
\item
go up by $z$;
\item
go along the second edge of $q_0$;
\item
go down by $z$;
\item
go along $p_2$
\end{itemize}
and so on.

In the middle of the process we will pass through $p_n$ (since $n$ is even).
Now we have already visited all vertices of $\Gamma_n$ with zero 3rd
coordinate. We are now at the point $(n,n,0)$. The process continues as
follows:

\begin{itemize}
\item
we do up by $z$;
\item
go along $q_n^{-1}$;
\item
go down by $z^{-1}$;
\item
go along the first edge of $p_0^{-1}$;
\item
go up by $z$;
\item
go along $q_{n-1}$
\end{itemize}
and so on. At the final step we will go along $q_0^{-1}$ and then go down
by $z$. We went home (to the initial point) visiting also all vertices of
$\Gamma_n$ with the 3rd coordinate equal to $1$.

Let us find the length of $r$.

For each $0\le i\le n$, we went along $p_i^{\pm1}$ exactly once. This took
$Nn$ steps; at each of them we went along an edge labelled by $v^{\pm1}$.
Since $|v|=4$, we used length $4Nn$ for the above intervals. Similarly, we
``spent" length $Nn$ for the intervals of the form $q_i^{\pm1}$. (Each
edge of $q_i$ is labelled by the word $u$ of length $1$.) So we have
$5Nn$ as a part of the whole length.

It is easy to verify that the number of times when we went up by $z$ is
equal to $2n+1$. The points on the plane with zero 3rd coordinate from
which we went up by $z$ form the following list: $(n,0)$, $(0,1)$, $(n,2)$,
$(0,3)$, \dots, $(0,n-1)$; $(n,n)$; $(n-1,0)$, $(n-2,n)$, $(n-3,0)$,
$(n-4,n)$, \dots, $(0,n)$. Also we went down by $z$ exactly the same
number of times. So for this part of the length we have $2(2n+1)|\xi|$.

But we also had a few single edges between going up and down (or down and
up). There were exactly $n$ edges of this kind labelled by $u^{\pm1}$ and
the same number of edges labelled by $v^{\pm1}$. So we need to add
$n(|u|+|v|)=5n$ to the length. Finally we have:

\be{est}
|r|=5Nn+5n+(4n+2)|\xi|=5(N^2-1)+(4N-2)|\xi|<5N^2+4N|\xi|.
\ee
Now for the ratio we have the following inequalities:

\be{ratio}
|r|/\Card(S_n)<(5N^2+4N|\xi|)/(2N^2)=2.5+2|\xi|/N<\lambda
\ee
provided $N>2|\xi|/(\lambda-2.5)$.

The proof is complete.
\endproof

{\bf Remark.}\ If we change the generating set of $F$ by $x_0$, $x_1$,
$x_2$, then it is possible to prove that the corresponding Cayley
graph does not satisfy $TS(\lambda)$ for $\lambda>1.5$. Indeed, in
this case one can take $u=x_1x_0^{-1}$, $v=x_2$. The proof goes
without any essential changes. The sum $|u|+|v|$ equals $3$ instead
of $5$ so $5/2$ replaces by $3/2$.

Another natural generating set of $F$ is $x_0$, $x_1$, $x_1x_0^{-1}$
(the third generator is a function, which is an image of $x_1$ under
conjugation by the transformation $t\mapsto1-t$ of the unit interval).
In this case, taking $u=x_1x_0^{-1}$ and $v=x_0^{-1}x_1x_0$, the sum
of length of the words in the new generating set becomes equal to $4$.
So the corresponding Cayley graph will not satisfy property $TS(\lambda)$
provided $\lambda>2$.
\vspace{2ex}

The author thanks Mark Sapir for helpful remarks.

\end{document}